\documentclass[a4paper,10pt,reqno]{amsart}
\usepackage[a4paper,tmargin=40mm,bmargin=35mm,hmargin=30mm]{geometry}
\usepackage[T1]{fontenc}
\usepackage{lmodern}
\usepackage{amsmath}
\usepackage{amssymb}
\usepackage{amsthm}
\usepackage{stmaryrd}
\usepackage{tikz}
\usepackage{booktabs}
\usepackage{multirow}
\usepackage{multicol}
\usepackage{here}
\usepackage{aliascnt}
\usepackage{hyperref}
\usepackage[noabbrev]{cleveref}

\newcommand{\NewTheorem}[2]{
	\newaliascnt{#1}{TheoremEnvironment}
	\newtheorem{#1}[#1]{#1}
	\aliascntresetthe{#1}
	\crefname{#1}{#1}{#2}
	\Crefname{#1}{#1}{#2}
}

\theoremstyle{definition}

\NewTheorem{Definition}{Definitions}
\NewTheorem{Remark}{Remarks}
\NewTheorem{Axiom}{Axioms}
\NewTheorem{Example}{Examples}
\NewTheorem{Observation}{Observations}
\NewTheorem{Convention}{Conventions}
\NewTheorem{Notation}{Notations}
\NewTheorem{Setting}{Settings}
\NewTheorem{Question}{Questions}
\NewTheorem{Answer}{Answers}
\NewTheorem{Conjecture}{Conjectures}
\NewTheorem{Problem}{Problems}
\NewTheorem{Solution}{Solutions}
\NewTheorem{Goal}{Goals}
\NewTheorem{Comment}{Comments}
\NewTheorem{Aim}{Aims}
\NewTheorem{Caution}{Cautions}
\NewTheorem{Exercise}{Exercises}

\theoremstyle{plain}
\NewTheorem{Proposition}{Propositions}
\NewTheorem{Lemma}{Lemmas}
\NewTheorem{Theorem}{Theorems}
\NewTheorem{Corollary}{Corollaries}

\crefname{enumi}{}{}
\Crefname{enumi}{}{}
\creflabelformat{enumi}{(#2#1#3)}
\crefname{enumii}{}{}
\Crefname{enumii}{}{}
\creflabelformat{enumii}{(#2#1#3)}
\crefname{enumiii}{}{}
\Crefname{enumiii}{}{}
\creflabelformat{enumiii}{(#2#1#3)}

\makeatletter
\renewcommand{\p@enumii}{}
\renewcommand{\p@enumiii}{}
\makeatother

\numberwithin{equation}{section}
\crefname{equation}{}{}
\Crefname{equation}{}{}
\creflabelformat{equation}{(#2#1#3)}

\newcommand{\SwapSymbols}[1]{
	\expandafter\let\expandafter\temporarysymbol\csname #1\endcsname
	\expandafter\let\csname #1\expandafter\endcsname\csname var#1\endcsname
	\expandafter\let\csname var#1\endcsname\temporarysymbol
}

\SwapSymbols{epsilon}
\SwapSymbols{phi}
\SwapSymbols{Gamma}
\SwapSymbols{Delta}
\SwapSymbols{Theta}
\SwapSymbols{Lambda}
\SwapSymbols{Xi}
\SwapSymbols{Pi}
\SwapSymbols{Sigma}
\SwapSymbols{Upsilon}
\SwapSymbols{Phi}
\SwapSymbols{Psi}
\SwapSymbols{Omega}

\newcommand{\cS}{\mathcal{S}}

\newcommand{\Modr}{\mathrm{Mod}\text{-}}

\newcommand{\To}{\longrightarrow}

\DeclareMathOperator{\Hom}{Hom}

\DeclareMathOperator{\Ext}{Ext}
\DeclareMathOperator{\Tor}{Tor}

\DeclareMathOperator{\Ann}{Ann}

\DeclareMathOperator{\Ker}{Ker}
\DeclareMathOperator{\Coker}{Coker}
\let\Im\relax
\DeclareMathOperator{\Im}{Im}

\DeclareMathOperator{\Supp}{Supp}

\DeclareMathOperator{\depth}{depth}

\title{A criterion for cofiniteness of modules}
\subjclass[2010]{13D45, 13E05}
\keywords{local cohomology, cofinite module}

\author{Mohammad Khazaei and Reza Sazeedeh}
\address{Department of Mathematics, Urmia University, P.O.Box: 165, Urmia, Iran}
\email{mohmmad.khazaei@gmail.com}
\address{Department of Mathematics, Urmia University, P.O.Box: 165, Urmia, Iran}
\email{rsazeedeh@ipm.ir}

\begin{document}

\begin{abstract}
Let $A$ be a commutative noetherian ring, $\frak a$ be an ideal of $A$, $m,n$ be non-negative integers and let $M$ be an $A$-module such that $\Ext^i_A(A/\frak a,M)$ is finitely generated for all $i\leq m+n$. We define a class $\cS_n(\frak a)$ of modules and we assume that $H_{\frak a}^s(M)\in\cS_{n}(\frak a)$ for all $s\leq m$. We show that $H_{\frak a}^s(M)$ is $\frak a$-cofinite for all $s\leq m$ if either $n=1$ or $n\geq 2$ and $\Ext_A^{i}(A/\frak a,H_{\frak a}^{t+s-i}(M))$ is finitely generated for all $1\leq t\leq n-1$, $i\leq t-1$ and $s\leq m$. If $A$ is a ring of dimension $d$ and $M\in\cS_n(\frak a)$ for any ideal $\frak a$ of dimension $\leq d-1$, then we prove that  $M\in\cS_n(\frak a)$ for any ideal $\frak a$ of $A$.
  	\end{abstract}

\maketitle
\tableofcontents

\section{Introduction}
Throughout this paper $A$ is a commutative noetherian ring, $\frak a$ is an ideal of $A$, $m,n$ are non-negative integer numbers and $M$ is an $A$-module, unless otherwise stated. Grothendieck \cite{G}, expos$\acute{e}$ 13, 1.2 conjectured that if $M$ is a finitely generated $A$-module, then $\Hom_A(A/\frak a,H_{\frak a}^i(M))$ is finitely generated, where $H_{\frak a}^i(M)$ is the $i$-th local cohomology of $M$ with respect to the ideal $\frak a$. 

The concept of cofiniteness of modules was defined for the first time by Hartshorne \cite{H}, giving a negative answer to the Grothendieck's conjecture and later was studied by the numerous authors \cite{BNS, DM, HK, MV, M1, M2, M3}. 

An $A$-module $M$ is said to be $\frak a$-{\it cofinite} if $\Supp M\subseteq V(\frak a)$ and $\Ext_A^i(A/\frak a,M)$ is finitely generated for all integers $i\geq 0$.
In \cite{NS}, the authors introduced a criterion for the cofiniteness of modules. In a roughly speaking, the criterion estimates how much a module is close to be cofinite. We recall that an $A$-module $M$ satisfies the condition P$_n(\frak a)$ if the following implication holds:
\begin{center}
P$_n(\frak a):\hspace{0.3cm}$ {\it If $\Ext_A^i(A/\frak a, M)$ is
finite for all $i\leq n$ and $\Supp(M)\subseteq V(\frak a)$,
then $M$ is $\frak a$-cofinite.}
 \end{center}
We denotes by $\cS_n(\frak a)$, the class of all modules satisfying the condition P$_n(\frak a)$. 

We start this paper by describing and getting some explanations of the class $\cS_n(\frak a)$, especially in the case where $n\leq 2$. We give several examples and results at the beginning of the paper. But, our main aim of this paper is to study the cofiniteness of local cohomology modules when they belong to $\cS_n(\frak a)$ without any condition on the ideal $\frak a$ or on the ring $A$. The first result is devoted for the case $n=1$ as follows.
\begin{Theorem}$\rm{(\cref{s1})}$
If $\Ext_A^i(A/\frak a,M)$ is finitely generated for all $i\leq m+1$ and $H_{\frak a}^i(M)\in\cS_1(\frak a)$ for all $i\leq m$, then $H_{\frak a}^i(M)$ is $\frak a$-cofinite for all $i\leq m$.
\end{Theorem}
In the case where $\dim A/\frak a=1$, Melkersson \cite{M3} showed that $\cS_1(\frak a)=\Modr A$ and so the theorem implies that $H_{\frak a}^i(M)$ is $\frak a$-cofinite for all $i$ whenever $\Ext^i(A/\frak a,M)$ is finitely generated for all $i\leq n+1$ where $n=\dim M$.  

Moreover, we have the following result for the class $\cS_n(\frak a)$ where $n\geq 2$.

\begin{Theorem}$\rm{(\cref{tt})}$
Assume that $m$ is a non-negative integer such that $\Ext^i(A/\frak a,M)$ is finitely generated for all $i\leq m+n$ and $H_{\frak a}^s(M)\in\cS_{n}(\frak a)$ for all $s\leq m$. If $\Ext_A^{i}(A/\frak a,H_{\frak a}^{t+s-i}(M))$ is finitely generated for all $1\leq t\leq n-1$, $i\leq t-1$ and $s\leq m$, then $H_{\frak a}^s(M)$ is $\frak a$-cofinite for all $s\leq m$. 
\end{Theorem}

As an application of this theorem, let $\Ext_A^i(A/\frak a,M)$ be finitely generated for all $i\leq m+2$ and let $H_{\frak a}^{i}(M)\in\cS_2(\frak a)$ for all $i\leq m$. We show that if $\Hom_A(A/\frak a,H_{\frak a}^i(M))$ is finitely generated for all $i\leq m+1$, then $H_{\frak a}^{i}(M)$ is $\frak a$-cofinite for all $i\leq m$.
 In the case where $\dim A/\frak a=2$ where $A$ is a local ring, Bahmanpour et al. \cite{BNS} showed that $\cS_2(\frak a)=\Modr A$ and so this application generalizes \cite[Theorem 3.7]{BNS} and \cite[Theorem 3.7]{NS}. For the case of $\dim A/\frak a=3$, assume that $\depth(\Ann M,A/\frak a)>0$ and $\Ext_A^i(A/\frak a,M)$ is finitely generated for all $i\leq 2$. Then we show that $\Gamma_{\frak a}(M)$ is $\frak a$-cofinite if $\Hom_A(A/\frak a, H_{\frak a}^i(M))$ is finitely generated for $i=0,1$. 
 
We give a result about those modules which their local cohomology modules are nonzero only in two consecutive numbers. To be more precise, assume that $t$ is a non-negative integer such that $\Ext^i_A(A/\frak a,M)$ is finitely generated for all $i\leq n+t+1$ and $H_{\frak a}^i(M)=0$ for all $i\neq t,t+1$. Then we show that $H_{\frak a}^{t+1}(M)\in\cS_n(\frak a)$ if and only if $H_{\frak a}^t(M)\in\cS_{n+2}(\frak a)$.

One of the substantial results in the local cohomology theory and cofiniteness is the change of ring principle. We show that this result holds for $\cS_n(\frak a)$ as well. Let $B$ be a finitely generated $A$-algebra and let $M$ be a $B$-module. Then we show that $M\in\cS_n(\frak a)$ if and only if $M\in\cS_n(\frak aB).$  

We end this paper by the following theorem which is a generalization of \cite[Theorem 2.3]{NS}:

\begin{Theorem}$\rm{(\cref{quot})}$
Let $A$ be a ring of dimension $d$ and $M\in\cS_n(\frak a)$ for any ideal $\frak a$ of dimension $\leq d-1$ (i.e $\dim A/\frak a\leq d-1$). Then $M\in\cS_n(\frak a)$ for any ideal $\frak a$ of $A$. 
\end{Theorem}

\section{The main results} 

We start this section by defining a class of $A$-modules which has an essential role in this paper. For the basic properties of local cohomology modules, we refer the readers to \cite{BS}.
\begin{Definition}
Let $n$ be a non-negative integer
and let $M$ be an $A$-module. We say that $M$ satisfies the {\it
condition P$_n(\frak a)$} if the following implication holds:
\begin{center}
P$_n(\frak a):\hspace{0.3cm}$ {\it If $\Ext_A^i(A/\frak a, M)$ is
finitely generated for all $i\leq n$ and $\Supp(M)\subseteq V(\frak a)$,\\
then $M$ is $\frak a$-cofinite.}
 \end{center}
We define a class of $A$-modules as follows 
$$\cS_n(\frak a)=\{M\in\Modr A| \hspace{0.1cm} M\hspace{0.1cm} \textnormal{satisfies the condition}\hspace{0.1cm} P_n(\frak a)\}.$$
 We observe that $\cS_0(\frak a)\subseteq \cS_1(\frak a)\subseteq \dots$. We also say that $A$ {\it satisfies the condition
P$_n(\frak a)$} if $\cS_n(\frak a)=\Modr A$ where $\Modr A$ denotes the category of $A$-modules. 
\end{Definition}

\medskip

In the rest of this section, we assume that $\frak a$ is an ideal of $A$, $n$ is a non-negative integer and $M$ is an $A$-module, unless otherwise stated. In order to describe the class $\cS_n(\frak a)$, we give several examples. The first example shows that the top local cohomology modules lie in $\cS_0(\frak a)$.

\medskip
\begin{Example}
Assume that $\frak a$ is an arbitrary ideal of $A$ and $M$ is an $A$-module of dimension $d$ where $\dim M$ means the dimension of $\Supp M$. Then $H_{\frak a}^d(M)$ is in $\cS_0(\frak a)$. To be more precise, if $\Hom_A(A/\frak a,H_{\frak a}^d(M))$ is a finitely generated $A$-module, then it follows from \cite[Theorem 3.11]{NS} that $H_{\frak a}^d(M)$ is artinian and so, since $\Hom_A(A/\frak a,H_{\frak a}^d(M))$ has finite length, according to \cite[Proposition 4.1]{M2}, the module $H_{\frak a}^d(M)$ is $\frak a$-cofinite.  
\end{Example}

The following example specifies some modules in $\cS_1(\frak a)$.
\medskip

\begin{Example}
Given an arbitrary ideal $\frak a$ of $A$, by virtue of \cite[Lemma 2.2]{BNS}, $M\in\cS_1(\frak a)$ for all modules $M$ with $\dim M\leq 1$. Especially, if $\dim A/\frak a=1$, then it follows from \cite[Theorem 2.3]{M3} that $\cS_1(\frak a)=\Modr A$. Furthermore, if $\dim A=2$, then it follows from \cite[Corollary 2.4]{NS} that $\cS_1(\frak a)=\Modr A$ for any ideal $\frak a$ of $A$. 
\end{Example}
 
For the class $\cS_2(\frak a)$, we have the following example. 
\medskip
\begin{Example}
Let $\frak a$ be an ideal of a local ring $A$ with $\dim A/\frak a=2$. It follows from \cite[Theorem 3.5]{BNS} that $\cS_2(\frak a)=\Modr A$. Furthermore, if $A$ is a local ring with $\dim A=3$, then it follows from \cite[Corollary 2.5]{NS} that $\cS_2(\frak a)=\Modr A$ for any ideal $\frak a$ of $A$. 
\end{Example}

In the above example we may have $\cS_0(\frak a)\neq\Modr A$ or $\cS_1(\frak a)\neq\Modr A$. 
\medskip

\begin{Example}
The Hartshorne's example \cite{H} shows that if $\frak a$ is an ideal of $A$ with $\dim A/\frak a=2$, then we may have both $\cS_0(\frak a)\neq\Modr A$ and $\cS_1(\frak a)\neq\Modr A$. More precisely, assume that $A=k[[x,y,u,v]]$ where $x,y,u,v$ are variables and $k$ is a field. let $\frak p=(x,u)$ and $M=A/(xy-uv)$. Then $H_{\frak p}^1(M)\notin\cS_1(\frak a)$ and $H_{\frak p}^2(M)\notin\cS_0(\frak a)$. To be more precise, $H_{\frak p}^i(A)=0$ for all $i\neq 2$ as $\depth(\frak p,A)=2$ and since $\depth(\frak p,M)=1$, we have $\Gamma_{\frak p}(M)=0$; and hence we have an exact sequence of modules 
$0\To H_{\frak p}^1(M)\To H_{\frak p}^2(A)\stackrel{xy-uv.}\To H_{\frak p}^2(A)\To H_{\frak p}^2(M)\To 0$. By \cite{H}, the module $H_{\frak p}^2(M)$ is not $\frak p$-cofinite and it follows from \cite[Proposition 2.5]{MV} that  $H_{\frak p}^2(A)$ is $\frak p$-cofinite, and hence $H_{\frak p}^1(M)$ is not $\frak p$-cofinite. We observe that $\Ext^i_A(A/\frak p,H_{\frak p}^1(M))$ is finitely generated for $i=0,1$ so that $H_{\frak p}^1(M)\notin\cS_1(\frak a)$. Moreover, $\Hom_A(A/\frak p,H_{\frak p}^2(M))$ is finitely generated so that $H_{\frak p}^2(M)\notin\cS_0(\frak a)$.   
\end{Example}

For the case $\dim A/\frak a=3$, we have the following result.

\medskip
\begin{Proposition}
Let $A$ be a local ring with $\dim A/\frak a=3$, $\depth(\Ann M,A/\frak a)>0$ and let $\Ext_A^i(A/\frak a,M)$ is finitely generated for all $i\leq 2$. If $\Hom_A(A/\frak a, H_{\frak a}^i(M))$ is finitely generated for $i=0,1$, then $\Gamma_{\frak a}(M)$ is $\frak a$-cofinite. 
\end{Proposition}
\begin{proof}
There exists an element $x\in\Ann M$ such that $x$ is an $A/\frak a$-sequence; and hence $\dim \frac{A}{xA+\frak a}=2$. It follows from \cite[Proposition 1]{DM} that $\Ext_A^i(A/xA+\frak a,M)$ is finitely generated for all $i\leq 2$. On the other hand, using \cite[ Proposition 2]{DM}, for each $i\geq 0$, the module $H_{\frak a}^i(M)$ is $\frak a$-cofinite if and only if $H_{xA+\frak a}^i(M)$ is $xA+\frak a$-cofinite. Set $\frak b=xA+\frak a$ and $\overline{M}=M/\Gamma_{\frak a}(M)$. Then there is an exact sequence of modules $0\to \overline{M}\To E\To N\To 0$ such that $E$ is an injective $A$-module with $\Gamma_{\frak a}(E)=0$ and so $\Gamma_{\frak b}(E)=0$ as $\frak a\subseteq \frak b$. By the assumption, $\Hom_A(A/\frak a,\Gamma_{\frak a}(N))\cong\Hom_A(A/\frak a,H_{\frak a}^1(M))$ is finitely generated and since $\frak a\subseteq \frak b$ and $\Gamma_{\frak b}(N)\subseteq\Gamma_{\frak a}(N)$, the module $\Hom_A(A/\frak b,\Gamma_{\frak b}(N))$ is finitely generated. Since $x\in\Ann M$, we have $\Gamma_{\frak a}(M)=\Gamma_{\frak b}(M)$; and hence $\Hom_A(A/{\frak b},\Gamma_{\frak b}(M))$ is finitely generated. Therefore, the isomorphisms $H_{\frak b}^1(M)\cong H_{\frak b}^1(\overline{M})\cong \Gamma_{\frak b}(N)$ imply that $\Hom_A(A/\frak b,H_{\frak b}^1(M))\cong\Hom_A(A/\frak b,\Gamma_{\frak b}(N))$ is finitely generated. Now, it follows from \cite[Theorem 3.7]{NS} that $\Gamma_{\frak b}(M)$ is $\frak b$-cofinite so that the first argument deduces that $\Gamma_{\frak a}(M)$ is $\frak a$-cofinite.  
\end{proof}

The following result establishes a relation between the classes $\cS_n(\frak a)$ and $\cS_n(\frak p_i)$ where $\frak p_i$ are the minimal prime ideals of $\frak a$ for $1\leq i\leq t$.
\medskip
\begin{Proposition}
Let $\frak p_1,\dots,\frak p_t$ be the minimal prime ideals of $\frak a$ and let $\Supp (M)\subseteq V(\frak p_1+\dots+\frak p_t)$. If $M\in\cS_n(\frak p_i)$ for each $1\leq i\leq t$, then $M\in\cS_n(\frak a)$. 
\end{Proposition}
\begin{proof}
Clearly $\Supp M\subseteq V(\frak p_i)\subseteq V(\frak a)$. Now assume that $\Ext_A^j(A/\frak a,M)$ is finitely generated for all $1\leq j\leq n$. It follows from \cite[Proposition 1]{DM} that $\Ext_A^j(A/\frak p_i,M)$ is finitely generated for all $1\leq j\leq n$ and all $1\leq i\leq t$; and hence using the assumption,  $\Ext_A^j(A/\frak p_i,M)$ is finitely generated for all $j\geq 0$ and all $1\leq i\leq t$. Then it follows from \cite[Corollary 1]{DM} that  $\Ext_A^j(A/\frak a,M)$ is finitely generated for all $j\geq 0$.
\end{proof}

\medskip

\begin{Proposition}\label{co}
Let $x\in\frak a$ and $\Supp M\subseteq V(\frak a)$ such that $(0:_Mx), M/xM\in\cS_1(\frak a)$. Then $M\in\cS_2(\frak a)$.  
\end{Proposition}
\begin{proof}
Assume that $\Ext_A^i(A/\frak a,M)$ is finitely generated for all $i\leq 2$. Applying the functor $\Hom_A(A/\frak a,-)$ to the exact sequences of modules $0\To (0:_Mx)\To M\To xM\To 0$ and $0\To xM\To M\To M/xM\To 0$, it is straightforward to see that 
$\Ext_A^i(A/\frak a,(0:_Mx))$ is finitely generated for $i=0,1$ and since $(0:_Mx)\in\cS_1(\frak a)$, we conclude that $(0:_Mx)$ is $\frak a$-cofinite. But this implies that $\Ext_A^i(A/\frak a, M/xM)$ is finitely generated for $i=0,1$ and since $M/xM\in\cS_1(\frak a)$, we conclude that 
$M/xM$ is $\frak a$-cofinite. It now follows from \cite[Corollary 3.4]{M2} that $M$ is $\frak a$-cofinite. 
\end{proof}
\medskip
For the local cohomology modules of a finitely generated $A$-module of dimension 3 we have the following result.
\begin{Proposition}
If $M$ is a finitely generated $A$-module of dimension $3$ such that $H_{\frak a}^2(M)\in\cS_0(\frak a)$, then $H_{\frak a}^1(M)\in\cS_2(\frak a)$.
\end{Proposition}
\begin{proof}
Assume that $\Ext^i(A/\frak a,H_{\frak a}^1(M))$ is finitely generated for $i\leq 2$. We may assume that $\Gamma_{\frak a}(M)=0$ and so $\frak a$ contains a non-zerodivisor $x$ of $M$. Application the functor $\Gamma_{\frak a}(-)$ to the exact sequence $0\To M\stackrel{x.}\To M\To M/xM\To 0$ gives the following exact sequence $$0\To\Gamma_{\frak a}(M/xM)\To H_{\frak a}^1(M)\stackrel{x.}\To H_{\frak a}^1(M)\To H_{\frak a}^1(M/xM)\To H_{\frak a}^2(M)\stackrel{x.}\To\dots.$$
Since $\Ext^i(A/\frak a,H_{\frak a}^1(M))$ is finitely generated for $i\leq 2$, it is straightforward to see that $\Hom_A(A/\frak a, H_{\frak a}^2(M))$ is finitely generated and so the assumption implies that $H_{\frak a}^2(M)$ is $\frak a$-cofinite. On the other hand $\dim M/xM=2$ and so it follows from \cite[Proposition 5.1]{M2} that $H_{\frak a}^2(M/xM)$ and $H_{\frak a}^3(M)$ are $\frak a$-cofinite and by virtue of \cite[Proposition 2.5]{MV}, the module $H_{\frak a}^1(M/xM)$ is $\frak a$-cofinite. Thus, it is straightforward to show that $H_{\frak a}^1(M)/xH_{\frak a}^1(M)$ is $\frak a$-cofinite and hence it follows from \cite[Corollary 3.4]{M2} that $H_{\frak a}^1(M)$ is $\frak a$-cofinite.
\end{proof}
\medskip

The following result is the first main theorem about cofiniteness of local cohomology modules when they lie in $\cS_1(\frak a)$.

\begin{Theorem}\label{s1}
If $\Ext_A^i(A/\frak a,M)$ is finitely generated for all $i\leq m+1$ and $H_{\frak a}^i(M)\in\cS_1(\frak a)$ for all $i\leq m$, then $H_{\frak a}^i(M)$ is $\frak a$-cofinite for all $i\leq m$.
\end{Theorem}
\begin{proof}
We proceed by induction on $m$. If $m=0$, then the isomorphism $\Hom_A(A/\frak a,\Gamma_{\frak a}(M))\cong \Hom_A(A/\frak a,M)$ and the exact sequence $0\To\Ext_A^1(A/\frak a, \Gamma_{\frak a}(M))\To \Ext_A^1(A/\frak a,M)$ imply that $\Ext_A^i(A/\frak a,\Gamma_{\frak a}(M))$ is finite for $i\leq 1$ and since $\Gamma_{\frak a}(M)\in\cS_1(\frak a)$, we deduce that $\Gamma_{\frak a}(M)$ is $\frak a$-cofinite. Now, suppose that $m>0$ and the result has been proved for all values $<m$. Considering $\overline{M}=M/\Gamma_{\frak a}(M)$, there is an exact sequence of modules $0\To \overline{M}\To E\To N\To 0$ such that $E$ is injective and $\Gamma_{\frak a}(E)=0$.  The case $m=0$ implies that $\Gamma_{\frak a}(M)$ is $\frak a$-cofinite so that $\Ext_A^i(A/\frak a,\overline{M})$ is finitely generated for all $i\leq m+1$. Thus the isomorphisms $\Ext^i_A(A/\frak a, N)\cong\Ext_A^{i+1}(A/\frak a,\overline{M})$ for all $i\geq 0$ implies that $\Ext^i_A(A/\frak a, N)$ is finitely generated for all $i\leq m$; furthermore $H_{\frak a}^i(N)\cong H_{\frak a}^{i+1}(M)\in\cS_1(\frak a)$ for all $i\leq m-1$. Now, the induction hypothesis implies that $H_{\frak a}^i(N)$ is $\frak a$-cofinite for all $i\leq m-1$ and the isomorphism $H_{\frak a}^i(N)\cong H_{\frak a}^{i+1}(M)$ for all $i\geq 0$, forces that $H_{\frak a}^i(M)$ is $\frak a$-cofinite for all $i\leq m$.  
\end{proof}

We now extend the above theorem for the class $\cS_n(\frak a)$ where $n\geq 2$.  
\medskip
\begin{Theorem}\label{tt}
Assume that $m$ is a non-negative integer such that $\Ext^i(A/\frak a,M)$ is finitely generated for all $i\leq m+n$ and $H_{\frak a}^s(M)\in\cS_{n}(\frak a)$ for all $s\leq m$. If $\Ext_A^{i}(A/\frak a,H_{\frak a}^{t+s-i}(M))$ is finitely generated for all $1\leq t\leq n-1$, $i\leq t-1$ and $s\leq m$, then $H_{\frak a}^s(M)$ is $\frak a$-cofinite for all $s\leq m$. 
\end{Theorem}
\begin{proof}
We proceed by induction on $m$. Assume that $m=0$ and consider the exact sequences $$0\To \Gamma_{\frak a}(M)\To M\To \overline{M}\To 0\hspace{0.2cm} (1);$$
		$$0\To\overline{M}\To E_0\To M_1\To 0\hspace{0.2cm} (1');$$
		$$.$$
		$$.$$
		$$.$$
	$$0\To \Gamma_{\frak a}(M_{i-1})\To M_{i-1}\To \overline{M_{i-1}}\To 0\hspace{0.2cm} (i);$$
		$$0\To\overline{M_{i-1}}\To E_{i-1}\To M_i\To 0\hspace{0.2cm} (i').$$
		In view of the exact sequence $(1)$ and the fact that $\Gamma_{\frak a}(M)\in\cS_{n}(\frak a)$, it suffices to show that $\Ext_A^i(A/\frak a,\overline{M})$ is finitely generated for all $i\leq n-1$. Fix $i\leq n-1$. The case $i=0$ is clear. The case  $i=1$, in view of $(1')$, the module $\Ext^1(A/\frak a,\overline{M})\cong\Hom_A(A/\frak a,M_1)\cong\Hom_A(A/\frak a,\Gamma_{\frak a}(M_1))\cong\Hom_A(A/\frak a,H_{\frak a}^1(M))$ is finitely generated by the assumption (consider $t=1$ and $i=0$). For $1<i\leq n-1$, using $(1')$ we have an isomorphism $\Ext_A^i(A/\frak a,\overline{M})\cong\Ext^{i-1}_A(A/\frak a,M_1)$. Now, using $(2)$, it suffices to show that $\Ext^{i-1}_A(A/\frak a,\Gamma_{\frak a}(M_1))$ and $\Ext^{i-1}_A(A/\frak a,\overline{M_1})$ are finitely generated. The module $\Ext^{i-1}_A(A/\frak a,\Gamma_{\frak a}(M_1))\cong \Ext^{i-1}_A(A/\frak a,H_{\frak a}^1(M))$ is finitely generated by the assumption (replace $t$ by $i-1$ and $i$ by $i-2$). Continuing this way, we have to prove that $\Ext_A^1(A/\frak a,\Gamma_{\frak a}(M_{i-1}))$ and $\Ext_A^1(A/\frak a,\overline{M_{i-1}})$ are finitely generated. Using the above exact sequence $\Ext_A^1(A/\frak a,\Gamma_{\frak a}(M_{i-1}))\cong \Ext_A^1(A/\frak a, H_{\frak a}^{i-1}(M))$ and so using the assumption, it is finitely generated (replace $t$ by $i$ and $i$ by $1$). On the other hand, in view of $(i')$, $\Ext_A^1(A/\frak a,\overline{M_{i-1}})\cong\Hom_A(A/\frak a,M_i)\cong \Hom_A(A/\frak a,\Gamma_{\frak a}(M_i))\cong \Hom_A(A/\frak a,H_{\frak a}^i(M))$ is finitely generated by the assumption (replace $t$ by $i$ and $i$ by $0$). Suppose that $m>0$ and the result has been proved for all values $<m$. By the induction hypothesis, $\Gamma_{\frak a}(M)$ is $\frak a$-cofinite so that $\Ext_A^i(A/\frak a,\overline{M})$ is finitely generated for all $i\leq m+n$ and so in view of $(1')$, the module $\Ext_A^i(A/\frak a, M_1)\cong\Ext_A^{i+1}(A/\frak a,\overline{M})$ is finitely generated for all $i\leq m+n-1$. On the other hand, for all $1\leq t\leq n-1$, $i\leq t-1$ and $s\leq m$, the modules $\Ext_A^{i}(A/\frak a,H_{\frak a}^{t+s-1-i}(M_1))\cong \Ext_A^{i}(A/\frak a,H_{\frak a}^{t+s-i}(M))$ are finitely generated and hence using the induction hypothesis $H_{\frak a}^s(M)\cong H_{\frak a}^{s-1}(M_1)$ is $\frak a$ cofinite for all $s\leq m$,
		\end{proof}
The following corollary is a generalization of \cite[Theorem 3.7]{NS} and \cite[Theorem 3.7]{BNS} without any conditions on $A$ and $\frak a$.
\medskip
\begin{Corollary}\label{coff}
Let $\Ext_A^i(A/\frak a,M)$ be finitely generated for all $i\leq m+2$ and let $H_{\frak a}^{i}(M)\in\cS_2(\frak a)$ for all $i\leq m$. If $\Hom_A(A/\frak a,H_{\frak a}^i(M))$ is finitely generated for all $i\leq m+1$, then $H_{\frak a}^{i}(M)$ is $\frak a$-cofinite for all $i\leq m$.
 \end{Corollary}
\begin{proof}
The proof is straightforward by the previous theorem considering $n=2$. 
\end{proof}

When the local cohomology modules of a module are nonzero only in two consecutive numbers, we have the following result.

\medskip

\begin{Proposition}
Let $t$ be a non-negative integer such that $H_{\frak a}^i(M)=0$ for all $i\neq t,t+1$ and let $\Ext^i_A(A/\frak a,M)$ be finitely generated for all $i\leq n+t+1$. Then $H_{\frak a}^{t+1}(M)\in\cS_n(\frak a)$ if and only if $H_{\frak a}^t(M)\in\cS_{n+2}(\frak a)$. 
\end{Proposition}
\begin{proof}
There is the Grothendieck spectral sequence $$E_2^{p,q}:=\Ext_A^p(A/\frak a, H_{\frak a}^q(M))\Rightarrow \Ext_A^{p+q}(A/\frak a,M).$$
For each $p$, consider the sequence $E_2^{p-2,t+1}\stackrel{d_2^{p-2,t+1}}\To E_2^{p,t}\stackrel{d_2^{p,t}}\To E_2^{p+2,t-1}$. By the assumption $E_2^{p+2,t-1}=0$ so that $E_3^{p,t}=\Ker d_2^{p,t}/\Im d_2^{p-2,t+1}=\Coker d_2^{p-2,t+1}$. Now consider the sequence 
$E_3^{p-3,t+2}\stackrel{d_3^{p-3,t+2}}\To E_3^{p,t}\stackrel{d_3^{p,t}}\To E_3^{p+3,t-2}$. Since $E_3^{p-3,t+2}$ and $E_3^{p+3,t-2}$ are the subquotients of $E_2^{p-3,t+2}$ and $ E_2^{p+3,t-2}$ respectively, the assumption implies that $E_3^{p-3,t+2}\stackrel{d_3^{p-3,t+2}}=E_3^{p+3,t-2}=0$ so that $E_4^{p,t}=E_3^{p,t}$. Continuing this manner, we deduce that $E_3^{p,t}=E_{\infty}^{p,t}$; and hence there is an exact sequence of modules $$E_2^{p-2,t+1}\To E_2^{p,t}\To E_{\infty}^{p,t}\To 0\hspace{0.2cm}(\dag).$$ Using a similar argument, we have another exact sequence of modules $$0\To E_{\infty}^{p,t+1}\To E_2^{p,t+1}\To E_2^{p+2,t}\hspace{0.2cm}(\ddag).$$ As $E_{\infty}^{p,t}$ and $E_{\infty}^{p,t+1}$ are the subquotients of $\Ext_A^{p+t}(A/\frak a,M)$ and $\Ext_A^{p+t+1}(A/\frak a,M)$ respectively, they are finitely generated for all $p\leq n$ by the assumption. Assume that
 $H_{\frak a}^{t+1}(M)\in\cS_n(\frak a)$ and $\Ext^p_A(A/\frak a, H_{\frak a}^t(M))$ is finitely generated for all $p\leq n+2.$ The exact sequence $(\ddag)$ implies that $\Ext_A^p(A/\frak a, H_{\frak a}^{t+1}(M))$ is finitely generated for all $p\leq n$; and hence $H_{\frak a}^{t+1}(M)$ is $\frak a$-cofinite. It now follows from $(\dag)$ that $H_{\frak a}^t(M)$ is $\frak a$-cofinite so that $H_{\frak a}^t(M)\in\cS_{n+2}(\frak a)$. The converse is obtained by a similar argument. 
\end{proof}
\medskip

\begin{Example}
Let $k$ be a field of characteristic $0$, and let $R=K[X_{ij}]$ for $1\leq i\leq 2$ and $1\leq j\leq 3$. Let $\frak p$ be the height two prime ideal generated by $2\times 2$ minors of the matrix $(X_{ij})$. As $\frak p$ is generated by $3$ elements and $A$ is domain, $H_{\frak p}^i(A)=0$ for all $i\neq 2,3$. Since $\Hom_A(A/\frak p, H_{\frak p}^3(A))$ is not finitely generated, we have $H_{\frak p}^3(A)\in\cS_0(\frak p)$ and hence the previous proposition implies that $H_{\frak p}^2(A)\in\cS_2(\frak p)$. 
\end{Example}

We show that the change of ring principle holds for $\cS_n(\frak a)$.

\medskip

\begin{Proposition}\label{cri}
Let $B$ be a finitely generated $A$-algebra and let $M$ be a $B$-module. Then $M\in\cS_n(\frak a)$ if and only if $M\in\cS_n(\frak aB).$  
\end{Proposition}
\begin{proof}
It is clear that $\Supp_A M\subseteq V(\frak a)$ if and only if $\Supp_B M\subseteq V(\frak aB)$.
Assume that $\Ext_B^i(B/\frak aB,M)$ is finitely generated for all $0\leq i\leq n$. Consider the Grothendieck spectral sequence $$E_2^{p,q}:=\Ext_B^p(\Tor_q^R(B,A/\frak a),M)\Rightarrow H^{p+q}=\Ext_A^{p+q}(A/\frak a,M).$$
By the assumption, $E_2^{p,0}$ is finitely generated for all $0\leq p\leq n$ and since $\Supp_B\Tor_q^A(B,A/\frak a)\subseteq V(\frak aB)$ for all $q\geq 0$, we deduce that $E_2^{p,q}$ is finitely generated for all $0\leq p\leq n$ and all $q\geq 0$ by \cite[Proposition 1]{DM}. For any $r>2$, the $B$-module $E_r^{p,q}$ is a subquotient of $E_{r-1}^{p,q}$ and so an easy induction yields that $E_r^{p,q}$ is finitely generated for all $r\geq 2$, $0\leq p\leq n$ and all $q\geq 0$ so that $E_{\infty}^{p,q}$ is finitely generated for all $0\leq p\leq n$ and all $q\geq 0$. For any $0\leq t\leq n$, there is a finite filtration 
$$0=\Phi^{t+1}H^t\subset \Phi^tH^t\subset\dots\subset\Phi^1H^t\subset \Phi^0H^t\subset H^t$$ 
such that $\Phi^pH^t/\Phi^{p+1}H^t\cong E_{\infty}^{p,t-p}$ where $0\leq p\leq t$. Since $E_{\infty}^{p,t-p}$ is finitely generated for all $0\leq p\leq t$ and $0\leq t\leq n$, we deduce that $H^t$ is finitely generated for all $0\leq t\leq n$ and since $M\in\cS_n(\frak a)$, we deduce that $M$ is $\frak a$-cofinite. Consequently, using \cite[Proposition 2]{DM}, the module $M$ is $\frak aB$-cofinite. Now, assume  that $M\in\cS_n(\frak aB)$ and $\Ext_A^i(A/\frak a,M)$ is finitely generated for all $0\leq i\leq n$. By induction on $i\leq n$, we show that $\Ext_B^i(B/\frak aB,M)$ is finitely generated $B$-module. For $i=0$, we have $\Hom_B(B/\frak aB,M)\cong \Hom_A(A/\frak a,M)$ is finitely generated. Now, assume that $i>0$ and the result has been proved for all values smaller than $i\leq n$. This means that $E_2^{p,0}=\Ext_B^p(B/\frak aB,M)$ is finitely generated for all $0\leq p<i$.  Since $\Supp_B\Tor_q^A(B,A/\frak a)\subseteq V(\frak aB)$, we conclude that $E_2^{p,q}$ is finitely generated for all $0\leq p<i$ and all $q$. The exact sequence $E_2^{i-2,1}\To E_2^{i,0}\To E_3^{i,0}\To 0$ and the induction hypothesis imply that $E_2^{i,0}$ is finitely generated if $E_3^{i,0}$ is finitely generated. Continuing this manner, we deduce that $E_2^{i,0}$ is finitely generated if $E_{\infty}^{i,0}$ is finitely generated. But there are the following filtration 
$$0=\Phi^{i+1}H^i\subset\dots\subset\Phi^1H^i\subset \Phi^0H^i\subset H^i$$  such that $E_{\infty}^{i,0}\cong \Phi^iH^i/\Phi^{i+1}H^i =\Phi^iH^i$ is a submodule of $H^i=\Ext_A^i(A/\frak a,M)$; and hence it is finitely generated. Therefore $\Ext_B^i(B/\frak aB,M)$ is finitely generated for all $0\leq i\leq n$ and since $M\in\cS_n(\frak aB)$, we deduce that $M$ is $\frak aB$-cofinite and it follows from \cite[Proposition 2]{DM} that $M$ is $\frak a$-cofinite.
\end{proof}

The following result is a generalization of \cite[Theorem 2.2]{NS}.
\medskip

\begin{Theorem}\label{quot}
Let $A$ be a ring of dimension $d$ and $M\in\cS_n(\frak a)$ for any ideal $\frak a$ of dimension $\leq d-1$ (i.e $\dim A/\frak a\leq d-1$). Then $M\in\cS_n(\frak a)$ for any ideal $\frak a$ of $A$. 
\end{Theorem}
\begin{proof}
Assume that $\frak a$ is an arbitrary ideal of $A$ such that $\Supp(M)\subseteq V(\frak a)$ and $\Ext_A^i(A/\frak a,M)$
is finitely generated for all $i\leq n$. We can choose a positive integer $t$ such that $(0:_A \frak
a^t)= \Gamma_{\frak a}(A)$. Put $\overline{A} = A/{\Gamma_{\frak
a}(A)}$ and $\overline{M} = M/{(0:_M \frak a^t)}$ which is an
$\overline{A}$-module. By \cite[Lemma 2.1]{NS}, the module $(0:_M\frak a^t)$ is finitely generated; and hence for any ideal $\frak b$ of $A$, it is clear that $M\in\cS_n(\frak b)$ if and only if $\overline{M}\in\cS_n(\frak b)$. Taking $\overline{\frak a}$ as the image of
$\frak a$ in $\overline{A}$, we have $\Gamma_{\overline{\frak
a}}({\overline{A}})=0$. Thus $\overline{\frak a}$ contains an
$\overline{A}$-regular element so that $\dim A/\frak a+\Gamma_{\frak
a}(A)={\rm dim}\overline{A}/\overline{\frak a}\leq d-1$. Thus the assumption implies that $M\in\cS_n(\frak a+\Gamma_{\frak
a}(A))$ and the previous arguments yields $\overline{M}\in\cS_n(\frak a+\Gamma_{\frak a}(A))$. Using the rings homomorphism $A\To \overline{A}$, it follows from \cref{cri} that $\overline{M}$ lies in $\cS_n(\frak a)$. In view of the exact sequence $0\To (0:_M\frak a^t)\To M\To \overline{M}\To 0$, the assumption on $M$ implies that $\Ext_A^i(A/\frak a,\overline{M})$ is finitely generated for all $i\leq n$; and hence $\overline{M}$ is $\frak a$-cofinite. Now the previous exact sequence implies that $M$ is $\frak a$-cofinite.  
\end{proof}


\providecommand{\bysame}{\leavevmode\hbox
to3em{\hrulefill}\thinspace}



\begin{thebibliography}{Saz21}

\bibitem {BNS} K. Bahmanpour, R. Naghipour and M. Sedghi, \emph{Cofiniteness with respect to ideals of small dimension},
  Algebr Represent Theor (2014), DOI 10.1007//s10468-014-9498-3.

\bibitem {BS} M. Brodmann, R.Y. Sharp, \emph{Local Cohomology: an Algebraic Introduction with Geometric Applications},
 Cambridge Univ. Press, Cambridge, UK (1998).
\bibitem {DM} D. Delfino, T. Marley,
\emph{Cofinite modules of local cohomology}, J. Algebra
\textbf{121(1)} (1997), 45--52.
\bibitem {G} A. Grothendieck,
\emph{Cohomologie locale des faisceaux coh$\acute{e}$rents et
th$\acute{e}$or$\grave{e}$mes de Lefschetz locaux et globaux (SGA
$2$)}, North-Holland, Amsterdam (1968).
\bibitem {H} R. Hartshorne,
\emph{Affine duality and cofiniteness}, Invent. Math. \textbf{9}
(1970), 145--164.
\bibitem {HK} C. Huneke, J. Koh, \emph{Cofiniteness and vanishing of local cohomology modules}, Math. Proc Camb. Phi. Soc, {\bf 110}(1991), 421-429.
\bibitem {MV} T. Marley, J. C. Vassilev, \emph{Cofiniteness and associated primes of local cohomology modules}, J. Algebra \textbf{256}(2002), 180-193.

\bibitem {M1}  L. Melkersson,
\emph{Properties of cofinite modules and applications to local
cohomology},
 Math. Proc. Camb. Phil. Soc. \textbf{125} (1999), 417--423.
\bibitem {M2}\label{2005} L. Melkersson, \emph{Modules cofinite with respect to an ideal},
J. Algebra \textbf{285} (2005), 649--668.
\bibitem {M3}\label{2012} L. Melkersson,
\emph{Cofiniteness with respect to ideals of dimension one}, J.
Algebra \textbf{372} (2012), 459--462.
\bibitem{NS} M. Nazari, R. Sazeedeh,
\emph{Cofiniteness with respect to two ideals and local cohomology}, Algebr Represent Theor, {\bf 22} (2019), 375-385.
\end{thebibliography}
\end{document}